\documentclass[11pt,reqno]{amsart}

\usepackage{lmodern}
\usepackage[T1]{fontenc}
\usepackage[utf8]{inputenc}
\usepackage[english]{babel}
\usepackage{microtype}

\usepackage{amsmath, amssymb, amsfonts, amsthm, mathtools}
\usepackage{physics}
\allowdisplaybreaks              % long align environments may split across pages
\numberwithin{equation}{section}

\usepackage{graphicx}
\graphicspath{{./images/}}
\usepackage{xcolor}
\usepackage{comment}

\theoremstyle{plain}
\newtheorem{theorem}{Theorem}[section]
\newtheorem{corollary}[theorem]{Corollary}
\newtheorem{lemma}[theorem]{Lemma}

\theoremstyle{definition}
\newtheorem{definition}[theorem]{Definition}

\theoremstyle{remark}

\usepackage[
  colorlinks=true,
  linkcolor=blue!60!black,
  citecolor=green!50!black,
  urlcolor=blue!60!black,
  breaklinks=true
]{hyperref}
\usepackage[capitalise,nameinlink]{cleveref}   

\title[Mixing Rates for Analytic Observables]{Mixing Rates for Analytic Observables: Toral Endomorphisms and Anosov Blaschke Deformations}

\author{Ekaterina Shchetka}
\address{Department of Mathematics, Rice University, Houston, TX 77005, USA}
\email{katia.shchetka@rice.edu}


\thanks{The author was partially supported by NSF grant DMS 2404309}
\makeatletter
\renewcommand{\MR}[1]{}
\makeatother

\makeatletter
\def\l@subsection{\@tocline{2}{0pt}{2pc}{5pc}{}}
\makeatother
\begin{document}

\begin{abstract}
We analyze exponential mixing for analytic observables and show that mixing at a double exponential rate occurs for genuinely nonlinear systems. Neither its presence nor its sharp exponent is stable under arbitrarily small area-preserving perturbations in the real-analytic topology.

We construct two families of area-preserving Anosov Blaschke product diffeomorphisms that converge in the real-analytic topology to the same hyperbolic toral automorphism. For nonzero parameters, both  deformations are topologically conjugate to the linear model but not $C^1$-conjugate. For all nonzero parameter values, the maps in the first family are mixing at a double exponential rate with the same sharp exponent, strictly smaller than that of the linear map. The sharp exponent is therefore not lower semicontinuous at the linear map. In the second family, mixing at a double exponential rate is lost altogether, although the maps remain exponentially mixing.

We also establish double exponential bounds for every ergodic toral endomorphism and uniform bounds of the same type for totally ergodic higher-rank toral actions. Our method combines the Fourier escape argument with Koopman operator estimates from a Hilbert space of analytic functions to its dual space of hyperfunctions.
\end{abstract}

\maketitle

\tableofcontents

%==============================================================================

\section{Introduction}

The decay of correlations for hyperbolic dynamical systems is studied through the spectral properties of composition and transfer operators.  In the real-analytic setting, the operator-theoretic study
of correlations and resonances goes back to Ruelle \cite{Rue76}. For real-analytic hyperbolic maps, correlation spectra and trace-class properties of composition
operators were considered by Rugh \cite{Rug92} and Faure--Roy
\cite{FR06}. Explicit resonance spectra for particular analytic
Anosov maps of $\mathbb T^2$ were computed in \cite{SBJ17,BJS22, PS23}, while
Adam \cite{A17} established the generic presence of nontrivial
resonances among real-analytic perturbations of hyperbolic toral
automorphisms.  More recently, J\'ez\'equel \cite{J24} obtained general
bounds on the distribution of resonances for real-analytic Anosov
diffeomorphisms.

For ergodic toral automorphisms, Lind \cite{Lin82} proved exponential mixing
for H\"older observables using Katznelson's Diophantine estimate
\cite{Kat71}.  He also observed that correlations
of smooth observables decay faster than every exponential; in this
sense, ergodic toral automorphisms are \textit{superexponentially mixing} on
$C^\infty(\mathbb T^d)$.  For \textit{analytic observables}, the exponential decay of Fourier coefficients
combines with the exponential escape of frequencies in Lind's argument
to produce \textit{double exponential}
decay in the number of iterates.  Estimates of this type for hyperbolic endomorphisms
and purely expanding endomorphisms were obtained by Brini--Siboni
\cite{BS01}. Related bounds for ergodic toral automorphisms occur in the work of Fannjiang--Wo{\l}owski \cite{FW03}, and of Feng--Iyer 
\cite{FI19} under the additional irreducibility assumption.

We call this phenomenon \emph{mixing at a double exponential rate}; see
Definition~\ref{d-dim:definition of mixing at a double exponential rate}. We show that every ergodic toral endomorphism is mixing at a double
exponential rate and obtain an explicit rate in terms of the least
expanding eigenvalue, the degree of the minimal polynomial on the
center-stable space, and the sizes of the relevant Jordan blocks (see Theorem~\ref{d-dim:ergodic auto is dem}). The result applies to reducible and non-diagonalizable
endomorphisms, including maps with eigenvalues on the unit circle and
non-unimodular maps.  The proof repeats Lind's Fourier escape argument, keeping
track of the rate through Katznelson's estimate. The same argument has higher-rank consequences.  In
Appendix~\ref{chpt:higher-rank-actions}, we obtain uniform double
exponential mixing estimates for totally ergodic $\mathbb Z^k$-actions
by toral automorphisms and for totally ergodic $\mathbb Z_+^k$-actions
by nonsingular toral endomorphisms.  The required uniform
Fourier escape estimates follow from
\cite{KK05,DK10,FKS13,L13}.

We next move from linear maps to their nonlinear Blaschke deformations. Our main result is that mixing at a double
exponential rate occurs for genuinely nonlinear systems, but neither its presence nor its \textit{sharp exponent} (see Definition~\ref{d-dim:definition of mixing at a double exponential rate}) is stable under arbitrarily small area-preserving perturbations in the real-analytic topology. This gives an invertible, area-preserving analogue on
$\mathbb T^2$ of the instability for uniformly expanding
finite Blaschke products on the circle \cite{Shc26}.

The main analytic idea is to view the Koopman precomposition operator as acting from the weighted Hilbert space of analytic functions (see Definition~\ref{d-dim:analytic space}) into its dual space of hyperfunctions (see Definition~\ref{d-dim:hyperfunction space}).

\begin{theorem}\label{intro:2d example}
Let $\mu$ denote the normalized Lebesgue measure on $\mathbb T^2$.
There exist a hyperbolic toral automorphism $T^{(0)}$ and two
one-parameter families
\[
T_\lambda^{(i)},\qquad i=1,2,\qquad \lambda\in(-1,1),
\]
of area-preserving real-analytic Anosov diffeomorphisms
with the following properties.

\begin{enumerate}
\item \textbf{\textup{Analytic convergence and nonconjugacy.}}

As $\lambda\to0$, the maps $T_\lambda^{(1)}$ and $T_\lambda^{(2)}$ converge to $T^{(0)}$ in the real-analytic topology. For every $\lambda\in(-1,1)$, both $T_\lambda^{(1)}$ and $T_\lambda^{(2)}$ are topologically conjugate to $T^{(0)}$. Nevertheless, for $\lambda\neq0$, neither $T_\lambda^{(1)}$
nor $T_\lambda^{(2)}$ is $C^1$-conjugate to $T^{(0)}$. 

\item \textbf{\textup{Mixing rates.}}
\begin{itemize}
\item
The linear map $T^{(0)}$ is mixing at a double exponential
rate with sharp exponent
\[
\frac12\log\Lambda_0,
\]
where $\Lambda_0>1$ is its expanding eigenvalue.

\item
For $\lambda\neq0$, the first deformation $T_\lambda^{(1)}$
is mixing at a double exponential rate with sharp exponent
\[
\frac12\log\Lambda_1,
\]
where $\Lambda_1$ is independent of $\lambda$ for $\lambda\neq0$ and $1<\Lambda_1<\Lambda_0$.

\item
For $\lambda\neq0$, the second deformation $T_\lambda^{(2)}$
is exponentially mixing but is not mixing at any
double exponential rate.
\end{itemize}
\item \textbf{\textup{Entropy.}}

For every $\lambda\in(-1,1)$ and $i=1,2$,
\[
h_{\mathrm{top}}(T_\lambda^{(i)})
=
h_{\mathrm{top}}(T^{(0)})
=
h_\mu(T^{(0)})
=
\log\Lambda_0.
\]
For $\lambda\neq0$, 
\[
\begin{aligned}
\log\Lambda_1
&<h_\mu(T_\lambda^{(1)})<\log\Lambda_0,\\
0
&<h_\mu(T_\lambda^{(2)})<\log\Lambda_0.
\end{aligned}
\]
Nevertheless,
\[
h_\mu(T_\lambda^{(i)})
\longrightarrow\log\Lambda_0
\qquad\text{as }\lambda\to0,
\qquad i=1,2.
\]
\end{enumerate}
\end{theorem}

The sharpness of the two double exponential rate exponents is
proved by constructing fixed analytic observables whose
correlations satisfy matching lower bounds along a
subsequence of iterates.

The nonlinear families from Theorem~\ref{intro:2d example} are two-dimensional Blaschke product
diffeomorphisms in the sense of Pujals--Shub \cite{PS08} and
Pujals--Roeder \cite{PR10}.  Their results provide the geometric setting and determine the common topological entropy. Both families remain in the homotopy class of $T^{(0)}$. 

Normalized Lebesgue measure is a measure of maximal entropy
for $T^{(0)}$, but not for either nonlinear deformation
when $\lambda\neq0$. These strict entropy inequalities
reflect classical entropy rigidity: for smooth
area-preserving Anosov diffeomorphisms of $\mathbb T^2$,
equality between the entropy of Lebesgue measure and
topological entropy is equivalent to smooth conjugacy
to the linear model (see \cite[Section~5b]{K07}). The first family shows that the sharp double exponential exponent drops from $\frac12\log\Lambda_0$ at $\lambda=0$
to $\frac12\log\Lambda_1$ for every $\lambda\neq0$.
It is therefore \textit{not lower semicontinuous}
at $T^{(0)}$, even in the real-analytic topology. In the second family, mixing at a double exponential rate is lost
altogether for every $\lambda\neq0$, although the maps
remain exponentially mixing. In both families, however, topological entropy remains constant and the entropy of Lebesgue measure is continuous at $\lambda=0$.
For general upper semicontinuity results for topological
entropy see Newhouse \cite{N89}.

Both Blaschke families fit within the spectral framework
of \cite{BJS22}, whose results compute the spectrum of the composition operator on suitable anisotropic spaces. For $\lambda\neq0$, their formulas show that $T_\lambda^{(1)}$ has no nontrivial resonances, whereas $T_\lambda^{(2)}$ does. Our estimates determine the sharp double exponential exponent for the first family. By Lemma~\ref{d-dim:analytic-to-smooth}, double exponential decay of correlations also implies that $T^{(0)}$ and $T_\lambda^{(1)}$ are superexponentially mixing on $C^\infty$-observables.

\section*{Acknowledgements}
The author is grateful to Alex Blumenthal, Jon Chaika, Jon DeWitt, David Fisher,  Andrey Gogolev, Boris Hasselblatt, Asaf Katz, Sarah Koch, Homin Lee, Roland Roeder, Julia Slipantschuk, Ralf Spatzier, Kurt Vinhage, and Amie Wilkinson for their insightful questions, helpful discussions, and comments on this work.

\section*{Use of AI tools}
GPT-5.6 Sol was used to check the grammar and language of the paper.
\section{Analytic Spaces and Fourier Escape Argument}
\label{section:analytic Fourier framework}

Throughout, $\mu$ denotes the normalized Lebesgue measure on $\mathbb T^d$.
Let \[T:\mathbb T^d\to\mathbb T^d\] preserve $\mu$.  Then $T$ is said to be \textit{(strongly) mixing} if, for any $f,g\in L_2(\mathbb{T}^d,\mu)$, the \textit{correlation
function}
\begin{equation*}\label{d-dim:correlation function}
\mathcal C_{f,g}(n)
:=\left|
\int_{\mathbb T^d}(f\circ T^n)\overline g\,d\mu
-\int_{\mathbb T^d}f\,d\mu
 \int_{\mathbb T^d}\overline g\,d\mu
\right|\to 0 \quad \text{as} \quad n\to\infty,
\end{equation*}
where $T^n$ denotes the $n$-fold composition of $T$ with itself. We study the \textit{quantitative mixing} for $T$ on the spaces of analytic functions $\mathcal{H}_\theta\subset L_2(\mathbb{T}^d,\mu)$.  

For $z=(z_1,\dots,z_d)\in\mathbb T^d$ and
$k=(k_1,\dots,k_d)\in\mathbb Z^d$, put
\[
 |k|:=|k|_1=\sum_{s=1}^d|k_s|,
 \qquad z^k:=z_1^{k_1}\cdots z_d^{k_d}.
\]

\begin{definition}[Analytic functions]
\label{d-dim:analytic space}
Let $0<\theta<1$.  The space of analytic functions $\mathcal H_\theta$
is defined as the completion of the trigonometric polynomials with respect
to
\[
 \|f\|_\theta^2
 :=\sum_{k\in\mathbb Z^d}|f_k|^2\theta^{-2|k|},
\]
where $f_k$ are the Fourier coefficients of $f$.
The closed subspace of elements with $f_0=0$ is denoted by
$\mathcal H_\theta^0$.
\end{definition}

\begin{definition}[Hyperfunctions]
\label{d-dim:hyperfunction space}
Let $0<\theta<1$.  The space of hyperfunctions
$\mathcal H_{\theta^{-1}}$ is defined as the completion of the
trigonometric polynomials with respect to
\[
 \|f\|_{\theta^{-1}}^2
 :=\sum_{k\in\mathbb Z^d}|f_k|^2\theta^{2|k|},
\]
where $f_k$ are the Fourier coefficients of $f$.
The closed subspace of elements with $f_0=0$ is denoted by
$\mathcal H_{\theta^{-1}}^0$.
\end{definition}

Under the pairing
\[
 \langle u,g\rangle=\sum_{k\in\mathbb Z^d}u_k\overline{g_k},
\]
we can identify $\mathcal H_{\theta^{-1}}$ isometrically with
$\mathcal H_\theta^*$.  In particular, the weighted Cauchy--Schwarz
inequality implies
\[
 |\langle u,g\rangle|
 \le \|u\|_{\theta^{-1}}\|g\|_\theta,
 \qquad
 u\in\mathcal H_{\theta^{-1}},\quad g\in\mathcal H_\theta.
\]
We also have the continuous inclusions
\[
 \mathcal H_\theta^0\hookrightarrow L_2^0(\mathbb T^d,\mu)
 \hookrightarrow\mathcal H_{\theta^{-1}}^0.
\]

\begin{definition}[Mixing at a double exponential rate]
\label{d-dim:definition of mixing at a double exponential rate}
Let $\gamma>0$. A map $T$ is said to be
\emph{mixing at a double exponential rate with exponent  $\gamma$} if, for some
$0<\theta<1$, there are $C,\beta>0$ such that for any $f,g\in\mathcal H_\theta$ and any $n\ge0$
\[
\mathcal C_{f,g}(n)
\le C e^{-\beta e^{\gamma n}}
\|f\|_\theta\|g\|_\theta.
\]
We say that $T$ is \textit{mixing at a double exponential rate} if this holds for some $\gamma>0$, and we call the exponent $\gamma$ \textit{sharp} if it cannot be replaced by any larger one. 
\end{definition}

If the estimate holds with exponent $\gamma$, it holds with every $0<\gamma'<\gamma$, so a sharp exponent, when it exists, is unique. By \cite[Lemma~B.1]{Shc26}, Definition~\ref{d-dim:definition of mixing at a double exponential rate} does not depend on $\theta$: if it holds for one, it holds for every $0<\theta<1$ with the same exponent $\gamma$. Although the lemma is stated for the circle, the proof applies verbatim to $\mathbb T^d$. 

The next lemma is a weighted analytic version of Lind's Fourier
argument \cite[Lemma~7.1]{Lin82}.  Its formulation in terms of
Fourier coefficients applies equally to nonlinear maps.

\begin{lemma}[Fourier escape argument]
\label{lem:Fourier escape argument}
Suppose that a sequence $R_n>0$ satisfies
\begin{equation}
\label{eq:abstract Fourier gap}
 \bigl((T^nz)^k\bigr)_j=0
 \quad\text{whenever}\quad
 0<|k|<R_n,\qquad 0<|j|<R_n.
\end{equation}
Then, for every $0<\theta<1$, there is $C_{\theta,d}>0$ such that
\begin{equation}
\label{eq:abstract Fourier gap correlation}
 \mathcal C_{f,g}(n)
 \le C_{\theta,d}\theta^{R_n/d}
 \|f\|_\theta\|g\|_\theta,
 \qquad f,g\in\mathcal H_\theta,\quad n\ge0.
\end{equation}
\end{lemma}

\begin{proof}
Define the \textit{Koopman} precomposition operator by
\[
 \mathcal Lf:=f\circ T, \quad f\in\mathcal{H}_\theta.
\]
Without loss of generality, assume that $f\in\mathcal{H}_\theta^0$. Then the weighted Cauchy--Schwarz inequality implies
\begin{align*}
 \mathcal C_{f,g}(n)
 &=|\langle\mathcal L^nf,g\rangle| \\
 &\le \|\mathcal{L}^nf\|_{\theta^{-1}}\|g\|_\theta \\
 &\le
 \|\mathcal L^n\|_{
 \mathcal H_\theta^0\to\mathcal H_{\theta^{-1}}^0}
 \|f\|_\theta\|g\|_\theta.
\end{align*}
The operator norm
is bounded by the Hilbert--Schmidt norm, and
\begin{align*}
 \|\mathcal L^n\|_{\mathcal H_\theta^0\to
      \mathcal H_{\theta^{-1}}^0}
 &\le \|\mathcal L^n\|_{\mathrm{HS}} \\
 &:=\sqrt{\sum_{\mathbf{0}\neq k\in\mathbb Z^d}
      \|\mathcal L^ne_k\|_{\theta^{-1}}^2} \\
 &=\sqrt{\sum_{\mathbf{0}\neq j,k\in\mathbb Z^d}
   \left|
   \left\langle\mathcal L^ne_k,e_j^*\right\rangle_
   {\mathcal H_{\theta^{-1}}}
   \right|^2},
\end{align*}
where
\[
 e_k:=\theta^{|k|}z^k,
 \qquad
 e_j^*:=\theta^{-|j|}z^j
\]
are orthonormal bases of $\mathcal H_\theta$ and
$\mathcal H_{\theta^{-1}}$, respectively. Next, we compute
\begin{align*}
 \left\langle\mathcal L^ne_k,e_j^*\right\rangle_
 {\mathcal H_{\theta^{-1}}}
 &=\theta^{|k|}\theta^{-|j|}\theta^{2|j|}
   \bigl((T^nz)^k\bigr)_j \notag\\
 &=\theta^{|j|+|k|}\bigl((T^nz)^k\bigr)_j.
\end{align*}
Thus, 
\begin{equation}\label{d-dim:main norm sum}
    \|\mathcal{L}^n\|_{\mathcal{H}_\theta^0\to \mathcal{H}_{\theta^{-1}}^0}\le \sqrt{\sum_{\mathbf{0}\neq j,k\in\mathbb{Z}^d}\theta^{2(|j|+|k|)}|((T^nz)^k)_j|^2} .
\end{equation}

 We decompose the sum in
\eqref{d-dim:main norm sum} into three parts:
\begin{equation}
\begin{split}
\|\mathcal{L}^n\|^2_{\mathcal{H}_\theta^0\to
\mathcal{H}_{\theta^{-1}}^0}
\leq
\bigg(&\sum_{0<|j|,|k|<R_n}
      +\sum_{\substack{|j|>0\\|k|\geq R_n}}
      +\sum_{\substack{|k|>0\\|j|\geq R_n}}\bigg)\theta^{2(|j|+|k|)}|((T^nz)^k)_j|^2.
\end{split}
\label{d-dim:three sums}
\end{equation}

Since $T^n(\mathbb T^d)\subseteq\mathbb T^d$,
\begin{equation}
\label{each coefficient is bounded by 1}
 \left|\bigl((T^nz)^k\bigr)_j\right|
 =\left|\int_{\mathbb T^d}(T^nz)^kz^{-j}\,d\mu(z)\right|
 \le1.
\end{equation}

By \eqref{eq:abstract Fourier gap}, the terms with
$0<|j|,|k|<R_n$ vanish. Hence
 \eqref{each coefficient is bounded by 1}
and Lemma~\ref{lem:higher-dimensional-geometric-series} with $\eta:=\theta^2$, applied to the
two remaining sums in \eqref{d-dim:three sums}, imply
\[
 \|\mathcal L^n\|_{\mathcal H_\theta^0\to
      \mathcal H_{\theta^{-1}}^0}
 \le C_{\theta,d}\theta^{R_n/d},
 \qquad
 C_{\theta,d}:=(2S_{\theta^2,d}D_{\theta^2,d})^{1/2}.
\]

This proves \eqref{eq:abstract Fourier gap correlation}.
\end{proof}

\section{Ergodic Toral Endomorphisms}
\label{ch3:auto}

Let $T\in\mathrm M(d,\mathbb Z)$ with $\det T\neq0$ be ergodic.  We use the same letter for the induced surjective endomorphism of $\mathbb T^d$: for
$z=(z_1,\dots,z_d)\in\mathbb T^d$,
\[
 Tz:=z^T
 =\bigl(z_1^{T_{11}}\cdots z_d^{T_{1d}},\dots,
         z_1^{T_{d1}}\cdots z_d^{T_{dd}}\bigr).
\]
Put $B:=T^{\mathsf t}$. Then for $k\in\mathbb Z^d$
\begin{equation}
\label{d-dim:character-action}
 (Tz)^k=z^{Bk},
 \qquad
 (T^nz)^k=z^{B^nk}.
\end{equation} Let
\[
 \mathbb R^d=E^{cs}\oplus E^u
\]
be the $B$-invariant splitting into generalized eigenspaces with
eigenvalues of modulus at most one and greater than one. For $v\in\mathbb R^d$, write
\[
 v=v_{cs}+v_u,
 \qquad v_{cs}\in E^{cs},\quad v_u\in E^u
\]
and according to the spectral
splitting, define the adapted norm
\[
 |v|_B:=|v_{cs}|+|v_u|.
\]
By equivalence of norms, there is $C_1=C_1(T)>0$ such that
\begin{equation}\label{equiv of norms}
 |v|\ge C_1|v|_B,
 \qquad v\in\mathbb R^d.
\end{equation}

Ergodicity implies that $E^u\neq\{0\}$; see
Appendix~\ref{chpt:katznelson-appendix}. Define
\[
 \Lambda:=
 \min\{ |\lambda|:\lambda\in\operatorname{spec}(B),\ |\lambda|>1\}.
\]
Let $p_{cs}$ be the minimal polynomial of $B|_{E^{cs}}$ and put
\begin{equation*}
 r:=\deg p_{cs}\le d-1,
\end{equation*}
with $p_{cs}=1$ and $r=0$ if $E^{cs}=\{0\}$.  If $s_\lambda$ denotes
the largest size of a Jordan block associated with $\lambda$, put
\begin{equation*}
 q:=\max_{\substack{\lambda\in\operatorname{spec}(B)\\
                     |\lambda|=\Lambda}}(s_\lambda-1).
\end{equation*}
\begin{theorem}
\label{d-dim:ergodic auto is dem}
Let $T\in\mathrm M(d,\mathbb Z)$, $\det T\neq0$, be ergodic.
For every $0<\theta<1$, there are $C,\beta>0$ such that
\begin{equation}
\label{d-dim:refined-correlation-bound}
\mathcal C_{f,g}(n)
\le C\exp\!\left(
-\beta n^{-q/(r+1)}\Lambda^{n/(r+1)}
\right)
\|f\|_\theta\|g\|_\theta
\end{equation}
for every $f,g\in\mathcal H_\theta$ and $n\ge1$.
In particular, for every
\[
0<\gamma<\frac{\log\Lambda}{r+1},
\]
the map $T$ is mixing at a double exponential rate $\gamma$. If $q=0$, then
\begin{equation}\label{the rate}
\gamma=\frac{\log\Lambda}{r+1}.
\end{equation}
\end{theorem}

Lind \cite{Lin82} proved that ergodic toral automorphisms are
exponentially mixing on H\"older observables using Katznelson's
Diophantine estimate \cite{Kat71}. We repeat his argument below, since we are interested in
the rate $\gamma$ explicitly given in terms of $\Lambda$ and $r$. He also observed
\cite[Section~7, Remark~(1)]{Lin82} that correlations of smooth
observables decay faster than any exponential.  Brini and Siboni
\cite{BS01} obtained double exponential decay for analytic observables
for hyperbolic and purely expanding toral maps.

Feng and Iyer \cite[Section~4]{FI19} obtain double exponential decay of
the energy of pulsed diffusion for ergodic
$A\in\mathrm{SL}_d(\mathbb Z)$ with irreducible characteristic
polynomial.  Their argument uses the simultaneous expansion of all
expanding conjugates, with exponent 
\[
\gamma=\frac{\log \left(\prod_{i=1}^d\max\{1,|\lambda_i|\}\right)}{d}.
\]
Within their irreducible class, if $E^u$ is one-dimensional, then $r+1=d$ and their exponent agrees with \eqref{the rate}.
When there are several unstable directions, \eqref{the rate} may be coarser because it retains only the weakest
unstable expansion. Its advantage is that it also covers reducible and
non-diagonalizable ergodic endomorphisms, including non-unimodular ones.

Fannjiang and Wo{\l}owski \cite{FW03} study noise-induced dissipation
for ergodic toral automorphisms. Their double exponential decay rates are expressed through the minimal dimensionally
averaged entropy $\widehat h$ rather than through $\Lambda$. 

\begin{lemma}[Lind's argument; \cite{Lin82}]
\label{d-dim:Lind's argument}
There exists $a=a(T)>0$ such that, with
\begin{equation}
\label{d-dim:refined-frequency-radius}
 R_n:=a n^{-q/(r+1)}\Lambda^{n/(r+1)},
\end{equation}
we have
\[
 \bigl((T^nz)^k\bigr)_j=0
 \quad\text{whenever}\quad
 0<|k|<R_n\quad\text{and}\quad 0<|j|<R_n.
\]
\end{lemma}

\begin{proof}
By \eqref{d-dim:character-action},
\begin{equation}
\label{linear map takes character to a character}
 \bigl((T^nz)^k\bigr)_j=\delta_{B^nk,j}.
\end{equation}
Write $k=k_{cs}+k_u$ according to the spectral splitting and use the
adapted norm
\[
 |k|_B:=|k_{cs}|+|k_u|.
\]

The equivalence of norms \eqref{equiv of norms}, unstable estimate
\eqref{d-dim:unstable-jordan-lower-bound}, and Katznelson's
Lemma~\ref{Katznelson} give constants $C_1,C_u,c_{\mathrm{K}}>0$ such that
\begin{align}
 |B^nk|
 &\ge C_1|B^nk_u|
 \ge C_1C_u n^{-q}\Lambda^n|k_u|\notag\\
 &\ge C_1C_u c_{\mathrm{K}} n^{-q}\Lambda^n|k|^{-r}.
\label{d-dim:refined-Bn-lower-bound}
\end{align}
Choose
\begin{equation*}
 a:=(C_1C_u c_{\mathrm{K}})^{1/(r+1)}.
\end{equation*}
If $0<|k|<R_n$, then
\eqref{d-dim:refined-Bn-lower-bound} gives $|B^nk|\ge R_n$.
Thus $B^nk\neq j$ whenever also $0<|j|<R_n$, and the conclusion follows
from \eqref{linear map takes character to a character}.
\end{proof}

\begin{proof}[Proof of Theorem~\ref{d-dim:ergodic auto is dem}]
Lemma~\ref{d-dim:Lind's argument} verifies the hypothesis of the
Fourier escape argument, Lemma~\ref{lem:Fourier escape argument}, with
$R_n$ given by \eqref{d-dim:refined-frequency-radius}.  Therefore
\eqref{d-dim:refined-correlation-bound} follows with
\[
 C:=C_{\theta,d},
 \qquad
 \beta:=-\frac{a}{d}\log\theta>0.
\]
The final conclusion follows by absorbing the polynomial factor into
$e^{\gamma n}$ for every
$\gamma<\log\Lambda/(r+1)$ and adjusting $\beta$.
\end{proof}

Next, we construct two nonlinear Blaschke deformations of a linear map.
One satisfies the Fourier escape hypothesis of
Lemma~\ref{lem:Fourier escape argument}, while the other does not. 
\section{Two-Dimensional Blaschke Deformations}\label{2d Blaschke}

We now pass from linear toral endomorphisms to a first nonlinear case.
We start with the hyperbolic matrix
\[
A_0=\begin{pmatrix}
5 & 2 \\
2 & 1 
\end{pmatrix}
\]
acting on a 2-torus \[\mathbb{T}^2=\{z=(z_1,z_2)\in\mathbb{C}^2: |z_1|=|z_2|=1\}\] by the formula 
\begin{align*}
T^{(0)}: \,  \mathbb{T}^2 &\to \mathbb{T}^2
\\
(z_1,z_2)&\mapsto(z_1^5z_2^2, z_1^2z_2).
\end{align*}
For a real parameter $\lambda\in(-1,1)$, we consider two deformations of
$T^{(0)}$ given by
\begin{align}
T_\lambda^{(i)}: \quad \mathbb{T}^2&\to \mathbb{T}^2, \qquad i=1,2, \nonumber
\\
T_\lambda^{(1)}: \quad (z_1,z_2)&\mapsto (z_1^3(B_\lambda(z_1))^2 z_2^2,\,z_1B_\lambda(z_1)z_2),\label{main-ex1}
\\ 
T_\lambda^{(2)}:\quad (z_1,z_2)&\mapsto (z_1(B_\lambda(z_1))^4 z_2^2,\,(B_\lambda(z_1))^2z_2) \label{main-ex}
\end{align}
where $B_\lambda$ is a Blaschke factor
\[
B_{\lambda}(z):=\frac{z-\lambda}{1-{\lambda}z}.
\]
When $\lambda=0$, the Blaschke factor reduces to $B_0(z)=z$, and hence
\[
T_0^{(i)}=T^{(0)}, \qquad i=1,2.
\]

For $\lambda\neq0$, the matrix of coordinatewise local degrees of
$T_\lambda^{(1)}$ at the fixed point $\mathbf 0\in\mathbb{C}^2$ is
\[
A_1:=\begin{pmatrix}3&2\\1&1\end{pmatrix}.
\]
Indeed,
\[
B_\lambda(z_1)=-\lambda+O(z_1),
\]
and hence
\[
T_\lambda^{(1)}(z_1,z_2)
=
\bigl(
(\lambda^2+O(z_1))z_1^3z_2^2,\,
(-\lambda+O(z_1))z_1z_2
\bigr).
\]
Since $\lambda\neq0$, the leading exponent vectors are $(3,2)$ and
$(1,1)$, which are the rows of $A_1$. 

\begin{theorem}\label{thm: 2d examples}
Let $T^{(0)}$ and $T_\lambda^{(i)}$, $i=1,2$, be defined as above. Then
\begin{enumerate}
\item For every $\lambda\in(-1,1)$, the maps $T^{(0)}$,
$T_\lambda^{(1)}$, and $T_\lambda^{(2)}$ are
Lebesgue measure preserving Anosov diffeomorphisms of $\mathbb{T}^2$.
\item
As $\lambda\to0$, both $T_\lambda^{(1)}$ and $T_\lambda^{(2)}$
converge to $T^{(0)}$ in the real-analytic topology.
\item For $\lambda\neq0$, neither $T_\lambda^{(1)}$ nor
$T_\lambda^{(2)}$ is $C^1$-conjugate to $T^{(0)}$.
\item For $\lambda\neq0$, the map $T_\lambda^{(2)}$ is exponentially mixing but not
mixing at a double exponential rate. 
\item  The map $T^{(0)}$ is mixing at the sharp double exponential rate
\[\frac12\log\Lambda_0,\] 
where $\Lambda_0=3+2\sqrt2$ is the expanding eigenvalue of $A_0$.
\item For $\lambda\neq0$, the map
$T_\lambda^{(1)}$ is mixing at the sharp double exponential rate
\[\frac12\log\Lambda_1,\] where 
$\Lambda_1=2+\sqrt3$ is the expanding eigenvalue of $A_1$.
\item The metric entropy of the linear map satisfies
\[
h_\mu\bigl(T^{(0)}\bigr)
=
h_{\mathrm{top}}\bigl(T^{(0)}\bigr)
=
\log\Lambda_0.
\]
For $\lambda\neq0$,
\[
\log\Lambda_1
<
h_\mu\bigl(T_\lambda^{(1)}\bigr)
<
h_{\mathrm{top}}\bigl(T_\lambda^{(1)}\bigr)
=
\log\Lambda_0,
\]
and
\[
0
<
h_\mu\bigl(T_\lambda^{(2)}\bigr)
<
h_{\mathrm{top}}\bigl(T_\lambda^{(2)}\bigr)
=
\log\Lambda_0.
\]
Nevertheless,
\[
h_\mu(T_\lambda^{(i)})
\longrightarrow\log\Lambda_0
\qquad\text{as }\lambda\to0,
\qquad i=1,2.
\]
\end{enumerate}
\end{theorem}

When no confusion can arise, we use $T^{(i)}$
instead of $T_\lambda^{(i)}$.

\subsection{Invariance of Lebesgue Measure and Anosov Property}
\begin{proof}[Proof of Lebesgue measure invariance in Part~1 of
Theorem~\ref{thm: 2d examples}]
Since $T^{(0)}$ is a toral automorphism and $\det(A_0)=1$, it preserves
$\mu$.  For the two deformations, the Fourier criterion in
Lemma~\ref{lm:Leb meas preserv} reduces the statement to proving that,
for every $k\in\mathbb Z^2$ and $i=1,2$,
\[
\int_{\mathbb T^2}(T^{(i)}(z_1,z_2))^k\,d\mu(z_1,z_2)
=\delta_{k,0}.
\]
For $T^{(1)}$, integration first in $z_2$ gives
\begin{align*}
&\int_{\mathbb T^2}
  (T^{(1)}(z_1,z_2))^k\,d\mu(z_1,z_2)\\
&\qquad=
\left(
\int_{\mathbb S^1}
z_1^{3k_1+k_2}B_\lambda(z_1)^{2k_1+k_2}\,d\mu(z_1)
\right)
\left(
\int_{\mathbb S^1}z_2^{2k_1+k_2}\,d\mu(z_2)
\right)\\
&\qquad=
\delta_{2k_1+k_2,0}
\int_{\mathbb S^1}
z_1^{3k_1+k_2}B_\lambda(z_1)^{2k_1+k_2}\,d\mu(z_1).
\end{align*}
Next we use that the condition $2k_1+k_2=0$ makes the power of $B_\lambda$ equal to
zero
\begin{align*}
\int_{\mathbb T^2}(T^{(1)}(z_1,z_2))^k\,d\mu(z_1,z_2)&=
\delta_{2k_1+k_2,0}
\int_{\mathbb S^1}z_1^{3k_1+k_2}\,d\mu(z_1)\\
&=
\delta_{2k_1+k_2,0}\delta_{3k_1+k_2,0}
\\&=\delta_{k,0}.
\end{align*}

  Similarly,
\begin{align*}
\int_{\mathbb T^2}(T^{(2)}(z_1,z_2))^k\,d\mu(z_1,z_2)&=
\int_{\mathbb S^1}
z_1^{k_1}B_\lambda(z_1)^{4k_1+2k_2}\,d\mu(z_1)
\int_{\mathbb S^1}z_2^{2k_1+k_2}\,d\mu(z_2)\\
&=
\delta_{2k_1+k_2,0}
\int_{\mathbb S^1}z_1^{k_1}\,d\mu(z_1)\\
&=
\delta_{2k_1+k_2,0}\delta_{k_1,0}
\\&=\delta_{k,0}.
\end{align*}
Thus both
nonlinear maps preserve $\mu$.
\end{proof}

\begin{proof}[Proof of Anosov property in Part 1 of Theorem~\ref{thm: 2d examples}]

All three maps are diffeomorphisms of $\mathbb T^2$ with their inverses given by the formulas 
\begin{align}\label{2d:inverse-maps}
\left(T^{(0)}\right)^{-1}(z_1,z_2)
&=(z_1z_2^{-2},z_1^{-2}z_2^5),\notag\\
\left(T^{(1)}\right)^{-1}(z_1,z_2)
&=\bigl(z_1z_2^{-2},
z_1^{-1}z_2^3(B_\lambda(z_1z_2^{-2}))^{-1}\bigr),\\ \notag
\left(T^{(2)}\right)^{-1}(z_1,z_2)
&=\bigl(z_1z_2^{-2},
z_2(B_\lambda(z_1z_2^{-2}))^{-2}\bigr).
\end{align}
The two nonlinear deformations are \textit{two-dimensional Blaschke products} in the sense
of Pujals and Shub \cite{PS08} and Pujals and Roeder \cite{PR10}, and the
inverse formulas \eqref{2d:inverse-maps} show that their restrictions to
$\mathbb T^2$ are \textit{Blaschke product diffeomorphisms}. 

Put $B:=B_\lambda(z_1)$ and
$B':=B_\lambda'(z_1)$.  The complex Jacobians are
\begin{equation}
\label{2d:complex Jacobians}
\begin{aligned}
DT^{(0)}(z)
&=\begin{pmatrix}
5z_1^4z_2^2&2z_1^5z_2\\
2z_1z_2&z_1^2
\end{pmatrix},\\
DT^{(1)}(z)
&=\begin{pmatrix}
z_1^2z_2^2(3B^2+2z_1BB')&2z_1^3B^2z_2\\
z_2(B+z_1B')&z_1B
\end{pmatrix},\\
DT^{(2)}(z)
&=\begin{pmatrix}
z_2^2(B^4+4z_1B^3B')&2z_1B^4z_2\\
2BB'z_2&B^2
\end{pmatrix}.
\end{aligned}
\end{equation}

For $|z_1|=|z_2|=1$, put $z_1=e^{ix_1}$ and $z_2=e^{ix_2}$, then the derivative in angular coordinates is
\[
 DT_x^{(i)}
 :=\operatorname{diag}\bigl(T^{(i)}(z)\bigr)^{-1}
 DT^{(i)}(z)\operatorname{diag}(z_1,z_2).
\]
The complex Jacobians in \eqref{2d:complex Jacobians} give
\begin{equation}
\label{2d:angular Jacobians}
 DT_x^{(1)}=
 \begin{pmatrix}
 3+2p(x_1)&2\\
 1+p(x_1)&1
 \end{pmatrix},
 \qquad
 DT_x^{(2)}=
 \begin{pmatrix}
 1+4p(x_1)&2\\
 2p(x_1)&1
\end{pmatrix},
\end{equation}
where \[
 p(x_1)
 :=\frac{z_1B'}{B}
 =\frac{1-\lambda^2}{|1-\lambda z_1|^2}.
\]
Since
\begin{equation}\label{one more useful bound}
0<\frac{1-|\lambda|}{1+|\lambda|}\le p(x)\le \frac{1+|\lambda|}{1-|\lambda|},
\end{equation}
all entries are positive, so the one-variable circle maps occurring in
the two coordinate formulas are orientation-preserving immersions.
Therefore, \cite[Theorem~3.2]{PS08} gives a dominated splitting for each
nonlinear map. An area-preserving surface diffeomorphism with a dominated splitting is uniformly hyperbolic \cite[p.~117]{BV02}. Hence $T^{(1)}$ and $T^{(2)}$ are Anosov.
The map $T^{(0)}$ is Anosov because it is induced by the hyperbolic
matrix $A_0$.
\end{proof}

\subsection{Real-analytic Perturbations of the Linear Map}

\begin{proof}[Proof of Part 2 of Theorem~\ref{thm: 2d examples}]
We will show that, for all sufficiently small $|\lambda|$,
both families admit holomorphic extensions to a fixed
complex neighborhood $U$ of $\mathbb T^2$ satisfying
\[
\sup_{z\in U}
\|T_\lambda^{(i)}(z)-T^{(0)}(z)\|
=O(|\lambda|),
\qquad i=1,2.
\]
This implies the real-analytic convergence.

Fix $R>1$ and set
\[
U=\left\{
(z_1,z_2)\in\mathbb C^2:
R^{-1}<|z_j|<R,\quad j=1,2
\right\}.
\]
Choose $0<\lambda_*<R^{-1}$.
For $|\lambda|\le\lambda_*$, the function $B_\lambda$
is holomorphic on a neighborhood of the closed disk
$|z|\le R$, and
\[
|B_\lambda(z)-z|
=
\left|
\frac{\lambda(z^2-1)}{1-\lambda z}
\right|
\le
\frac{R^2+1}{1-\lambda_*R}\,|\lambda|.
\]
Thus $B_\lambda(z)=z+O(|\lambda|)$ uniformly
on this disk. Substitution into the formulas \eqref{main-ex1} and \eqref{main-ex} gives
\begin{equation*}
\sup_{z\in U}
\|T_\lambda^{(i)}(z)-T^{(0)}(z)\|
=O(|\lambda|),
\qquad i=1,2.
\end{equation*}
\end{proof}

\subsection{Failure of \texorpdfstring{$C^1$}{C1}-conjugacy}

\begin{proof}[Proof of Part 3 of Theorem~\ref{thm: 2d examples}]

First, we compute that all three maps $T^{(i)}$, $i=0,1,2$, have the same fixed points on $\mathbb{T}^2$
\[
(1,1), \quad (1,-1),\quad (-1,1), \quad (-1,-1).
\]
Indeed, for $T^{(1)}$ the second fixed-point equation gives
$z_1B_\lambda(z_1)=1$, which is equivalent to $z_1^2=1$; the first
equation then gives $z_2^2=1$.  For $T^{(2)}$, the second fixed-point
equation gives $B_\lambda(z_1)^2=1$.  Since $B_\lambda$ is a bijection of
$\mathbb T$ and $B_\lambda(\pm1)=\pm1$, this again gives $z_1^2=1$, and
the first equation gives $z_2^2=1$.

To check if $T^{(i)}$, $i=1,2$, are $C^1$-conjugate to $T^{(0)}$, we need to compare the derivatives at corresponding fixed points.
Indeed, suppose that a $C^1$ diffeomorphism $h:\mathbb{T}^2\to\mathbb{T}^2$ satisfies
\[
h\circ T^{(i)}=T^{(0)}\circ h.
\]
If $p$ is a fixed point of $T^{(i)}$ and $q:=h(p)$, then $q$ is a fixed
point of $T^{(0)}$.  The chain rule gives
\[
Dh_p\circ D(T^{(i)})_p
=D(T^{(0)})_q\circ Dh_p.
\]
Since $Dh_p$ is invertible, the derivatives at $p$ and $q$ are
similar and therefore have the same eigenvalues.  This is the periodic-data obstruction to $C^1$ conjugacy; see, for example, \cite{HK95}. 

The matrices $DT$ computed below are the complex Jacobians.  At a fixed
point $z=(z_1,z_2)\in\mathbb{T}^2$, the derivative of the restriction to $\mathbb T^2$ in angular coordinates
is
\[
P^{-1}DT(z)P,
\qquad
P:=\operatorname{diag}(z_1,z_2).
\]
It is therefore similar to $DT(z)$ and has the same eigenvalues. A
conjugacy may permute the four fixed points, but the computation below shows
that the spectrum of $DT^{(0)}$ is the same at all four of them.

For $(z_1,z_2)=(\pm1,\pm1)$, we have $z_1^2=z_2^2=1$, so we compute
\[
DT^{(0)}(z_1,z_2)=\begin{pmatrix}
    5 & 2s\\
2s & 1 
\end{pmatrix},
\]
where $s=z_1z_2=\pm1$. The characteristic polynomial is independent of the sign $s$ and is equal to
\[
x^2-6x+1.
\] 
Hence, the eigenvalues of $DT^{(0)}$ are
\[
x_\pm=3\pm2\sqrt{2}.
\]

Next, for $z_1=\pm1$ and $B=B_\lambda$,
\[
B(z_1)=z_1,\qquad  z_1B(z_1)=1, \qquad B'(z_1)=\left(\frac{1+\lambda}{1-\lambda}\right)^{z_1}.
\]
Hence, for $(z_1,z_2)=(\pm1,\pm1)$,
\begin{align*}
DT^{(1)}(z_1,z_2)&=\begin{pmatrix}
    z_1^2z_2^2\left(3B^2+2z_1BB'\right) & 2z_1^3B^2z_2\\
z_2(B+z_1B') & z_1B
\end{pmatrix} \\&=\begin{pmatrix}
    3+2a & 2s\\
s\left(1+a\right) & 1
\end{pmatrix},
\end{align*}
where $s=z_1z_2$ and
\[
a=\left(\frac{1+\lambda}{1-\lambda}\right)^{z_1}\neq 1 \quad \text{ for }\quad \lambda\neq 0.
\]

The characteristic polynomial is independent of the sign $s$ and is equal
to
\[
x^2-\left(4+2a\right)x+1.
\]
Thus, the eigenvalues of $DT^{(1)}(z_1,z_2)$ are of the form
\begin{equation}\label{eigenvalues}
x_{\pm}=2+a\pm\sqrt{(2+a)^2-1}\neq 3\pm 2\sqrt{2} \quad \text{for}\quad \lambda\neq 0.
\end{equation}
For $T^{(2)}$, at the same fixed points we obtain
\[
DT^{(2)}(z_1,z_2)
=\begin{pmatrix}
1+4a&2s\\
2sa&1
\end{pmatrix}.
\]
Its characteristic polynomial is
\[
x^2-(2+4a)x+1.
\]
For $\lambda\neq0$, we have $a\neq1$ and the eigenvalues of $DT^{(2)}(z_1,z_2)$ satisfy
\begin{equation}\label{T2 ev}
x_\pm
=
1+2a
\pm\sqrt{(1+2a)^2-1}\neq 3\pm2\sqrt2.
\end{equation}
Consequently, neither $T^{(1)}$ nor $T^{(2)}$ is
$C^1$-conjugate to $T^{(0)}$ when $\lambda\neq0$.
\end{proof}

\subsection{Lower Bound on Exponential Mixing}

\begin{proof}[Proof of Part 4 of Theorem~\ref{thm: 2d examples}] 

By Part~1, $T_\lambda^{(2)}$ is a Lebesgue measure preserving Anosov
diffeomorphism homotopic to $T^{(0)}$. The Franks--Manning conjugacy
theorem \cite{F69,M74} implies that it is topologically
conjugate to $T^{(0)}$ and therefore topologically mixing. The
Bowen--Ruelle theory proves exponential decay of correlations for H\"older observables with respect to Lebesgue measure
\cite{Bow75}.  We now show that this decay cannot be uniformly
double exponential.

For $0<|\lambda|<1$, put
\begin{align*}
T:=T^{(2)}:\quad \mathbb{T}^2&\to\mathbb{T}^2\\
(z_1,z_2)&\mapsto (z_1(B_\lambda(z_1))^4 z_2^2,\,(B_\lambda(z_1))^2z_2).
\end{align*}

 We follow the same strategy as in Lemma~3.2 of \cite{Shc26}: we find two observables $f$ and $g$ such that the correlation decay is exactly exponential in $n$, thus proving that it cannot decay at a double exponential rate.

For $n\ge0$, the second coordinate of $T^nz$ has the form
\[
(T^nz)_2=\prod_{j=0}^{n-1} \left(B_\lambda\left((T^jz)_1\right)\right)^2\,z_2.
\]
We choose $f(z_1,z_2)=g(z_1,z_2)=z_2$. Both observables have zero mean. Compute the correlation function
\begin{align*}
\mathcal{C}_{f,g}(n)&=\left|\int_{\mathbb{T}^2} (f\circ T^n)\overline{g} \,d\mu\right|\\&=\left|\int_{\mathbb{T}^2} \prod_{j=0}^{n-1} \left(B_\lambda\left((T^jz)_1\right)\right)^2 z_2 \, z_2^{-1}\, d\mu(z)\right| \\&=
\left|\int_{\mathbb{T}^2} \prod_{j=0}^{n-1} \left(B_\lambda\left((T^jz)_1\right)\right)^2 \,d\mu(z_1)\,d\mu(z_2)\right|.
\end{align*}

For $n\ge0$, the integrand 
\[\prod_{j=0}^{n-1} \left(B_\lambda\left((T^jz)_1\right)\right)^2\] is holomorphic in a neighborhood of the closed bidisc $\{z\in\mathbb{C}^2: |z_1|\le1, |z_2|\le 1\}$. Thus, by applying multivariable Cauchy integral formula and using that $T(\mathbf{0})=\mathbf{0}$, we get 
\begin{align*}
\mathcal{C}_{f,g}(n)&=\left|\prod_{j=0}^{n-1} \left(B_\lambda\left((T^j\mathbf{0})_1\right)\right)^2\right|\\&
=\left|\prod_{j=0}^{n-1} \left(B_\lambda\left(0\right)\right)^2\right|\\&=\left|B_\lambda(0)\right|^{2n}\\&=|-\lambda|^{2n}\\&=|\lambda|^{2n}.
\end{align*} 
Since $|\lambda|^{2n}$ cannot satisfy an upper bound of
the form
\[
C\exp(-\beta e^{\gamma n})
\]
for any $C,\beta,\gamma>0$, the map $T_\lambda^{(2)}$ is not mixing at a
double exponential rate.
\end{proof}

\subsection{Sharpness of the Exponent for the Linear Map}

\begin{proof}[Proof of Part 5 of Theorem~\ref{thm: 2d examples}]

The upper bound follows from
Theorem~\ref{d-dim:ergodic auto is dem}.  Indeed, the expanding
eigenvalue of
\[
A_0=\begin{pmatrix}5&2\\2&1\end{pmatrix}
\]
is
\[
\Lambda_0=3+2\sqrt2.
\]
Moreover, the center-stable space is one-dimensional and $A_0$ is
diagonalizable, so $r=1$ and $q=0$.  Therefore, for every
$0<\theta<1$, there are $C,\beta>0$ such that
\[
\mathcal C_{f,g}(n)
\le
C\exp\left(-\beta\Lambda_0^{n/2}\right)
\|f\|_\theta\|g\|_\theta,
\qquad f,g\in\mathcal H_\theta.
\]
Thus $T^{(0)}$ is mixing at a double exponential rate
\[
\frac12\log\Lambda_0.
\]

We now prove that this rate is sharp. We follow the strategy of Theorem~5.1 in \cite{Shc26}. For $m\geq0$, define
\[
k_m:=A_0^{-m}e_2,
\qquad
j_m:=A_0^me_2,
\qquad
e_2:=\begin{pmatrix}0\\1\end{pmatrix}.
\]
Since the eigenvalues of $A_0$ are
$\Lambda_0$ and $\Lambda_0^{-1}$, there is a constant
$C_0>0$ such that
\[
|k_m|+|j_m|\leq C_0\Lambda_0^m.
\]

For a fixed $0<\theta<1$, define
\begin{align*}
f(z)
&:=\frac{\sqrt3}{2}
\sum_{m=0}^{\infty}2^{-m}\theta^{|k_m|}z^{k_m},\\
g(z)
&:=\frac{\sqrt3}{2}
\sum_{m=0}^{\infty}2^{-m}\theta^{|j_m|}z^{j_m}.
\end{align*}
Then $f,g\in\mathcal H_\theta^0$ and
\[
\|f\|_\theta=\|g\|_\theta=1.
\]
All their Fourier coefficients are nonnegative. 
For the linear map $T^{(0)}$, the correlation function equals
\[
\mathcal C_{f,g}(n)
=
\left|
\sum_{k\in\mathbb Z^2\setminus\{\mathbf0\}}
f_k\,\overline{g_{A_0^nk}}
\right|.
\]
Since
\[
A_0^{2n}k_n
=A_0^{2n}A_0^{-n}e_2
=A_0^ne_2
=j_n,
\]
the corresponding term in the correlation at time $2n$ gives
\begin{align*}
\mathcal C_{f,g}(2n)
&\geq
\frac34\,4^{-n}\theta^{|k_n|+|j_n|}\\
&\geq
\exp\left(-C\Lambda_0^n\right)\\
&=
\exp\left(-C\left(\sqrt{\Lambda_0}\right)^{\,2n}\right)
\end{align*} for some $C>0$. Therefore no double exponential rate
strictly larger than
\[
\frac{\log\Lambda_0}{2}
\]
can hold. 
\end{proof}

\subsection{The Nonlinear Analog of Lind's Argument}
The goal of this section is to prove for $T_\lambda^{(1)}$ a nonlinear
analog of Lind's argument from
Lemma~\ref{d-dim:Lind's argument}.

Throughout this subsection, fix $0<|\lambda|<1$ and put
\[
T:=T_\lambda^{(1)},\qquad
T(z_1,z_2)
=
\bigl(
z_1^3B_\lambda(z_1)^2z_2^2,\,
z_1B_\lambda(z_1)z_2
\bigr).
\]
Let
\begin{equation}\label{matrix B}
A:=A_1=
\begin{pmatrix}
3&2\\
1&1
\end{pmatrix},
\qquad
B:=A^{\mathsf t}=
\begin{pmatrix}
3&1\\
2&1
\end{pmatrix},
\end{equation}
and denote the expanding eigenvalue of $A$ (and $B$) by
\begin{equation}\label{def:Lambda}
\Lambda:=\Lambda_1=2+\sqrt3.
\end{equation}

For every $k\in\mathbb Z^2$,
\[
(Tz)^k
=
z^{Bk}B_\lambda(z_1)^{(Bk)_2}.
\]
Iterating this identity gives
\begin{equation}
\label{2d:iterate character formula}
(T^nz)^k
=
z^{B^nk}
\prod_{r=0}^{n-1}
B_\lambda\bigl((T^rz)_1\bigr)^{(B^{n-r}k)_2}.
\end{equation}

Define the symmetric unstable cone
\[
U:=
\left\{
k\in\mathbb Z^2\setminus\{\mathbf0\}:
(3k_1+k_2)(2k_1+k_2)\geq0
\right\}.
\]

\begin{lemma}
\label{2d:unstable cone escape}
There exists $a_u=a_u(A)>0$ such that, for every $n\geq1$,
$k\in U$, and $j\in\mathbb Z^2$ satisfying
\[
0<|j|<a_u\Lambda^n,
\]
we have
\[
\bigl((T^nz)^k\bigr)_j=0.
\]
\end{lemma}

\begin{proof}
First, suppose that
\begin{equation}\label{positive}
3k_1+k_2\geq0,
\qquad
2k_1+k_2\geq0.
\end{equation}
Then $Bk\in\mathbb Z_{\geq0}^2$.  Since $B$ has positive entries,
\[
B^mk\in\mathbb Z_{\geq0}^2,
\qquad m\geq1.
\]
In particular, all the exponents occurring in the product in
\eqref{2d:iterate character formula} are nonnegative. Thus the product in \eqref{2d:iterate character formula} is holomorphic in an open neighborhood of the bidisc. Cauchy's integral formula therefore gives
\[
\bigl((T^nz)^k\bigr)_j=0
\]
whenever
\[
j_1<(B^nk)_1
\qquad\text{or}\qquad
j_2<(B^nk)_2.
\]
Since both coordinates of $B^nk$ are nonnegative, this implies
\[
\bigl((T^nz)^k\bigr)_j=0
\qquad\text{whenever}\qquad
|j|<|B^nk|.
\]

The cone $U$ is separated from the stable eigendirection of the
hyperbolic matrix $B$.  Therefore, since $\Lambda$ is the expanding
eigenvalue of $B$, there is $a_u>0$ such that
\[
|B^nk|\geq a_u\Lambda^n|k|,
\qquad k\in U,\quad n\geq1.
\]
Since every $k\in U$ is a nonzero integer vector,
$|k|\geq1$. Hence
\[
|B^nk|\geq a_u\Lambda^n.
\]
Therefore,
\[
\bigl((T^nz)^k\bigr)_j=0
\qquad\text{whenever}\qquad
|j|<a_u\Lambda^n.
\]

It remains to consider the other component of $U$, on which
\[
3k_1+k_2\leq0,
\qquad
2k_1+k_2\leq0.
\]
Then $-k$ satisfies \eqref{positive}.  Since
\begin{equation}
\label{2d:conjugate Fourier coefficients}
\bigl((T^nz)^k\bigr)_j
=
\overline{\bigl((T^nz)^{-k}\bigr)_{-j}},
\end{equation}
applying the argument to $-k$ and $-j$ proves the result.
\end{proof}

Define the complementary stable cone
\[
S:=
\left\{
k\in\mathbb Z^2\setminus\{\mathbf0\}:
(3k_1+k_2)(2k_1+k_2)<0
\right\}.
\]

\begin{lemma}
\label{2d:stable cone enters unstable cone}
There exists $a_s=a_s(A)>0$ such that, for every $n\geq1$ and
$k\in S$ satisfying
\[
0<|k|<a_s\Lambda^n,
\]
we have
\[
\bigl((T^nz)^k\bigr)_j=0
\qquad\text{for every }j\in S.
\]
\end{lemma}

\begin{proof}
By \eqref{2d:conjugate Fourier coefficients} and the symmetry
$-S=S$, it is enough to consider the component of $S$ on which
\[
3k_1+k_2>0,
\qquad
2k_1+k_2<0.
\]
For such $k$, we have $k_1>0$ and $k_2<0$.

Let
\begin{equation}\label{matrix D}
\mathcal I(z_1,z_2):=(z_1,z_2^{-1}),
\qquad
D:=\begin{pmatrix}1&0\\0&-1\end{pmatrix}.
\end{equation}
Then
\[
(\mathcal Iz)^m=z^{Dm},
\qquad m\in\mathbb Z^2.
\]
Using the inverse formula \eqref{2d:inverse-maps} and
\[
B_\lambda(\zeta^{-1})=B_\lambda(\zeta)^{-1},
\]
we obtain
\[
\mathcal I\circ\widetilde T
:=
T^{-1}\circ\mathcal I,
\qquad
\widetilde T(z_1,z_2)
=
\bigl(
z_1z_2^2,\,
z_1z_2^3B_\lambda(z_1z_2^2)
\bigr).
\]
The map $\widetilde T$ extends holomorphically to the bidisc, maps the
bidisc into itself, and satisfies
\[
\widetilde T(\mathbf0)=\mathbf0.
\]

Since $T$ and $\mathcal I$ preserve Lebesgue measure, for every
$j\in\mathbb Z^2$ we have
\begin{align*}
\bigl((T^nz)^k\bigr)_j
&=
\left\langle (T^nz)^k,z^j\right\rangle_{L_2}\\
&=
\left\langle z^k,(T^{-n}z)^j\right\rangle_{L_2}
\\&=
\left\langle (\mathcal Iz)^k,(T^{-n}(\mathcal Iz))^j\right\rangle_{L_2}
\\
&=\left\langle z^{Dk},(\mathcal I(\widetilde{T}^nz))^j\right\rangle_{L_2}
\\&=
\left\langle z^{Dk},(\widetilde T^nz)^{Dj}\right\rangle_{L_2}
\\&=\overline{
\bigl((\widetilde T^nz)^{Dj}\bigr)_{Dk}
}.
\end{align*}

Suppose first that $j$ belongs to the component of $S$ opposite to the
one containing $k$.  Then
\[
Dk\in\mathbb Z_{>0}^2,
\qquad
-Dj\in\mathbb Z_{>0}^2.
\]
The holomorphic property of $\widetilde T$ therefore
gives
\[
\bigl((\widetilde T^nz)^{Dj}\bigr)_{Dk}=0.
\]

Now suppose that $j$ belongs to the same component of $S$ as $k$, then $Dk,Dj\in\mathbb Z_{>0}^2$. The map $\widetilde T$ satisfies the hypotheses of
Lemma~\ref{2d:unstable cone escape}, with
\[
\widetilde B
=DB^{-1}D
=\begin{pmatrix}1&1\\2&3\end{pmatrix}.
\]
Applying that lemma to the index $Dj\in\mathbb{Z}_{>0}^2$ shows that there exists $a_s>0$ such that
\[
\bigl((\widetilde T^nz)^{Dj}\bigr)_{Dk}=0 \qquad\text{whenever}\qquad 0<|Dk|<a_s\Lambda^n.
\]
Since $|Dk|=|k|$, this finishes the proof.
\end{proof}

\begin{lemma}[Nonlinear Lind's argument]
\label{2d:nonlinear Fourier escape}
There exists $a=a(T)>0$ such that, for every $n\geq1$ and
$j,k\in\mathbb Z^2$ satisfying
\[
0<|k|<a\Lambda^{n/2},
\qquad
0<|j|<a\Lambda^{n/2},
\]
we have
\[
\bigl((T^nz)^k\bigr)_j=0.
\]
\end{lemma}

\begin{proof}
Choose
\[
a:=\min\{a_s,a_u\Lambda^{-1/2},\Lambda^{-1/2}\}>0.
\]
For $n=1$, there are no $k\in\mathbb Z^2$ satisfying
\[
0<|k|<a\Lambda^{1/2}\leq1,
\]
so it is enough to consider $n\geq2$.

If $k\in U$, the conclusion follows directly from
Lemma~\ref{2d:unstable cone escape}.  Suppose, therefore, that
$k\in S$. We use first $\lceil n/2\rceil$ iterates to move $k\in S$ into the unstable cone $U$, and then to apply the unstable-cone escape estimate during the remaining
$\lfloor n/2\rfloor$ iterates. More precisely, if
\[
0<|k|<a\Lambda^{n/2},
\]
then
\[
0<|k|<a_s\Lambda^{\lceil n/2\rceil}.
\]
Lemma~\ref{2d:stable cone enters unstable cone} gives an
$L_2$-convergent expansion
\[
(T^{\lceil n/2\rceil}z)^k
=
\sum_{m\in U}c_{m,k,n}z^m.
\]
Applying the remaining $\lfloor n/2\rfloor$ iterates and taking the
$j$-th Fourier coefficient gives
\[
\bigl((T^nz)^k\bigr)_j
=
\sum_{m\in U}c_{m,k,n}
\left\langle
(T^{\lfloor n/2\rfloor}z)^m,z^j
\right\rangle_{L_2}.
\]
Moreover,
\[
0<|j|<a\Lambda^{n/2}
\leq a_u\Lambda^{\lfloor n/2\rfloor}.
\]
Every inner product in the sum therefore vanishes by
Lemma~\ref{2d:unstable cone escape}.  This proves the claim.
\end{proof}

\begin{proof}[Proof of the upper bound in Part 6 of Theorem~\ref{thm: 2d examples}]
Lemma~\ref{2d:nonlinear Fourier escape} verifies the hypothesis of the
Fourier escape argument, Lemma~\ref{lem:Fourier escape argument}, with
\[
R_n=a\Lambda_1^{n/2}.
\]
Consequently, for every $0<\theta<1$, there are $C,\beta>0$ such that
\[
\mathcal C_{f,g}(n)
\leq
C\exp\left(-\beta\Lambda_1^{n/2}\right)
\|f\|_\theta\|g\|_\theta,
\qquad
f,g\in\mathcal H_\theta.
\]
Thus $T_\lambda^{(1)}$ is mixing at a double exponential rate
\[
\frac12\log\Lambda_1.
\]
\end{proof}
Notice that this argument applies to every $0<|\lambda|<1$ and does not
use perturbation theory for small $\lambda$.

\subsection{Sharpness of the Exponent for Nonlinear Map} We now prove that the rate found in the previous section is sharp.  We follow the strategy of
\cite[Theorem~5.1]{Shc26}. We find $f,g\in\mathcal{H}_\theta^0$ and $C',\beta'>0$ such that 
\[
\mathcal{C}_{f,g}(n)\ge C'\exp\left(-\beta'\Lambda_1^{n/2}\right)\|f\|_\theta\|g\|_\theta
\]
for infinitely many $n$. 
Throughout this section $B$ is defined in \eqref{matrix B} and $\Lambda:=\Lambda_1$ is defined in \eqref{def:Lambda}.

For $m\geq0$, define
\[
k_m:=B^{-m}e_2,
\qquad
j_m:=B^me_2,
\qquad
e_2:=\begin{pmatrix}0\\1\end{pmatrix}.
\] 
Since the eigenvalues of $B$ are $\Lambda$ and $\Lambda^{-1}$, there is
$C_0>0$ such that
\begin{equation}\label{2d:nonlinear-balanced-frequencies}
|k_m|+|j_m|
\leq C_0\Lambda^m.
\end{equation}

\begin{lemma}
\label{lem:nonlinear-distinguished-coefficient}
For every $n\geq1$, define
\[
c_n:=
\bigl((T^{2n}z)^{k_n}\bigr)_{j_n}.
\]
Then there is $C>0$ such that
\begin{equation}
\label{2d:distinguished-coefficient-lower-bound}
|c_n|
\geq
e^{-C\Lambda^n}.
\end{equation}
Moreover, if $0\leq l,r\leq n$,
\begin{equation}
\label{2d:nonlinear-triangularity}
\bigl((T^{2n}z)^{k_l}\bigr)_{j_r}=0 \quad\text{unless}\quad l=r=n.
\end{equation}
\end{lemma}
\begin{proof}
We first note that
\begin{equation}
\label{2d:positive Blaschke exponents}
(B^se_2)_2>0,
\qquad s\in\mathbb Z.
\end{equation}
For $s\geq0$, this follows from the positivity of $B$. 
Let $D$ be defined in \eqref{matrix D}. Since
\[
\widetilde{B}=DB^{-1}D
=
\begin{pmatrix}
1&1\\
2&3
\end{pmatrix}
\]
has positive entries, for $s\geq1$ we have
\[
B^{-s}e_2
=
D(DB^{-1}D)^sDe_2
\]
also has positive second coordinate. This proves
\eqref{2d:positive Blaschke exponents}.

At time $2n$, apply \eqref{2d:iterate character formula} with
$k=k_n$
\[
(T^{2n}z)^{k_n}
=
z^{B^{2n}k_n}
\prod_{r=0}^{2n-1}
B_\lambda\bigl((T^rz)_1\bigr)^{(B^{2n-r}k_n)_2}.
\]
Since
\[
B^{2n}k_n=j_n,
\]
and all the Blaschke exponents are nonnegative by
\eqref{2d:positive Blaschke exponents}, the Cauchy formula gives
\begin{align*}
c_n
:&=
\bigl((T^{2n}z)^{k_n}\bigr)_{j_n}
\\&=\prod_{r=0}^{2n-1}
B_\lambda(0)^{(B^{2n-r}k_n)_2}\\&=
\prod_{r=0}^{2n-1}
B_\lambda(0)^{(B^{n-r}e_2)_2}\\&=
\prod_{s=1-n}^{n}
B_\lambda(0)^{(B^se_2)_2}.
\end{align*} 
Since
\[
B_\lambda(0)=-\lambda
\]
and
\[
\|B^s\|\leq C_1\Lambda^{|s|},
\]
there is $C>0$ such that
\[
|c_n|
\ge
e^{-C\Lambda^n}.
\]

Next we prove \eqref{2d:nonlinear-triangularity}.
Indeed, all the Blaschke exponents in
\eqref{2d:iterate character formula} with $k=k_l$ are nonnegative, so the Fourier
support of $(T^{2n}z)^{k_l}$ satisfies
\[
\operatorname{supp}\bigl((T^{2n}z)^{k_l}\bigr)
\subseteq
B^{2n}k_l+\mathbb Z_{\geq0}^2=B^{2n-l}e_2+\mathbb Z_{\geq0}^2.
\]
Since $B-I$ has nonnegative entries, the sequence
\[
B^qe_2,\qquad q\geq0,
\]
is coordinatewise nondecreasing and has no repeated terms.  
 Since
\[
2n-l\geq n\geq r,
\]
the frequency $j_r=B^re_2$ can occur only if
\[
2n-l=r,
\]
which implies $l=r=n$.  This proves
\eqref{2d:nonlinear-triangularity}.

\end{proof}

\begin{proof}[Proof of sharpness in Part 6 of Theorem~\ref{thm: 2d examples}]
Fix $0<\tau<1$ and define
\begin{align*}
f(z)
:=
\sum_{l=0}^{\infty}
\tau^{\Lambda^l}\theta^{|k_l|}z^{k_l} \quad \text{and}\quad
 g(z)
:=
\sum_{r=0}^{\infty}
\tau^{\Lambda^r}\theta^{|j_r|}z^{j_r}.
\end{align*}
We have $f,g\in\mathcal H_\theta^0$ and
\[
\|f\|_\theta
=
\|g\|_\theta
=
\sqrt{\sum_{r=0}^{\infty}\tau^{2\Lambda^r}}
<\infty.
\]

By Lemma~\ref{lem:nonlinear-distinguished-coefficient},
\begin{align}
\mathcal{C}_{f,g}(2n)&=\left|\left\langle
f\circ T^{2n}, g
\right\rangle_{L_2}\right|
\notag\\&=
\left|\sum_{l,r\geq0}
\tau^{\Lambda^l+\Lambda^r}
\theta^{|k_l|+|j_r|}
\bigl((T^{2n}z)^{k_l}\bigr)_{j_r}\right|
\notag\\
&=\left|
\tau^{2\Lambda^n}
\theta^{|k_n|+|j_n|}c_n
+
\sum_{\substack{l,r\geq0\\l>n\ \mathrm{or}\ r>n}}
\tau^{\Lambda^l+\Lambda^r}
\theta^{|k_l|+|j_r|}
\bigl((T^{2n}z)^{k_l}\bigr)_{j_r}\right|.
\label{2d:correlation-main-term-and-tail}
\end{align}

By \eqref{2d:nonlinear-balanced-frequencies} and
\eqref{2d:distinguished-coefficient-lower-bound}, there is
$C_{\theta,\lambda}>0$ such that 
\begin{equation*}
\theta^{|k_n|+|j_n|}|c_n|
\geq
e^{-C_{\theta,\lambda}\Lambda^n}.
\end{equation*}
Therefore, the main term in \eqref{2d:correlation-main-term-and-tail} satisfies
\[
\tau^{2\Lambda_1^n}
\theta^{|k_n|+|j_n|}|c_n|\ge
\bigl(\tau^2e^{-C_{\theta,\lambda}}\bigr)^{\Lambda^n}.
\]

Every Fourier coefficient of $(T^{2n}z)^{k_l}$ has modulus at most
one.  Therefore, the modulus of last sum in \eqref{2d:correlation-main-term-and-tail} is bounded above by
\begin{equation*}
\sum_{\substack{l,r\geq0 \\ l>n\ \mathrm{or}\ r>n}}
\tau^{\Lambda^l+\Lambda^r}
\leq
\frac{2\sum_{r=0}^{\infty}\tau^{\Lambda^r}}{1-\tau}\,
\tau^{\Lambda^{n+1}}.
\end{equation*}
Since $\Lambda>2$, we may choose $\tau>0$ sufficiently small that the tail is at most half of the main term.
In  other words,
\begin{align*}
\mathcal C_{f, g}(2n)
&\geq\frac12\bigl(\tau^2e^{-C_{\theta,\lambda}}\bigr)^{\Lambda^n}\\
&\geq
\exp\bigl(-\beta'\Lambda^n\bigr).
\end{align*}

After choosing $C'=\frac{1}{\|f\|_\theta\|g\|_\theta}$, we get for $m=2n\ge 0$
\[
\mathcal{C}_{f,g}(m)\ge C'\exp\left(-\beta'\Lambda^{m/2}\right)\|f\|_\theta\|g\|_\theta.
\]
Therefore no double exponential rate strictly larger than
\[
\frac12\log\Lambda
\]
can hold.
\end{proof}

\subsection{Metric and Topological Entropy}
Recall that $\mu$ denotes normalized Lebesgue measure on
$\mathbb T^2$.  We write $h_\mu(T)$ for the metric entropy of $T$
with respect to $\mu$.

\begin{proof}[Proof of Part~7 of Theorem~\ref{thm: 2d examples}]
By \cite[Proposition~A.1]{PR10}, all three maps have the same
topological entropy
\begin{equation}
\label{2d:common-topological-entropy}
h_{\mathrm{top}}\bigl(T_\lambda^{(i)}\bigr)
=h_{\mathrm{top}}\bigl(T^{(0)}\bigr)
=
\log\Lambda_0,
\qquad i=1,2.
\end{equation}

For the linear automorphism $T^{(0)}$, the positive Lyapunov exponent
is $\log\Lambda_0$.  Hence Pesin's entropy formula gives
\[
h_\mu\bigl(T^{(0)}\bigr)
=
h_{\mathrm{top}}\bigl(T^{(0)}\bigr)
=
\log\Lambda_0.
\]

We now consider the nonlinear maps.  By Part~1 of
Theorem~\ref{thm: 2d examples}, they are smooth
area-preserving Anosov diffeomorphisms. Pesin's entropy
formula \cite{P77} gives
\begin{equation}
\label{2d:Pesin-entropy-formula}
h_\mu\bigl(T_\lambda^{(i)}\bigr)
=
\int_{\mathbb T^2}\chi_i^+(x)\,d\mu(x),
\qquad i=1,2,
\end{equation}
where $\chi_i^+(x)$ denotes the positive Lyapunov exponent.

Let
\[
A_2:=
\begin{pmatrix}
1&2\\
0&1
\end{pmatrix}.
\]
Then the spectral radii
\[
\rho(A_1)=\Lambda_1,
\qquad
\rho(A_2)=1.
\]
By \eqref{one more useful bound},
\[
p_*:=\min_{x\in\mathbb T}p(x)>0.
\]
Put
\[
M_1:=A_1+p_*\begin{pmatrix}2&0\\1&0\end{pmatrix},
\qquad
M_2:=A_2+p_*\begin{pmatrix}4&0\\2&0\end{pmatrix}.
\]
Using \eqref{2d:angular Jacobians}, monotonicity of products of
nonnegative matrices, and Gelfand's formula for spectral radius, we obtain
\[
\chi_1^+(x)
\geq
\log\rho\left(
M_1
\right)
>
\log\rho(A_1)
=
\log\Lambda_1,
\]
and
\[
\chi_2^+(x)
\geq
\log\rho\left(
M_2
\right)
>
\log\rho(A_2)
=
0
\]
for almost every $x$.   
Thus
\eqref{2d:Pesin-entropy-formula} gives
\[
h_\mu\bigl(T_\lambda^{(1)}\bigr)>\log\Lambda_1,
\qquad
h_\mu\bigl(T_\lambda^{(2)}\bigr)>0.
\]

The variational principle and
\eqref{2d:common-topological-entropy} give
\[
h_\mu\bigl(T_\lambda^{(i)}\bigr)
\leq
h_{\mathrm{top}}\bigl(T_\lambda^{(i)}\bigr)
=
\log\Lambda_0.
\]

For $\lambda\neq0$, equality would imply, by entropy
rigidity for smooth area-preserving Anosov diffeomorphisms
of $\mathbb T^2$ \cite[Section~5b]{K07}, that
$T_\lambda^{(i)}$ is smoothly conjugate to its linear
model $T^{(0)}$. This contradicts Part~3. Thus
\[
h_\mu(T_\lambda^{(i)})
<h_{\mathrm{top}}(T_\lambda^{(i)})
=\log\Lambda_0,
\qquad i=1,2,\quad \lambda\neq0.
\]

Finally, as $\lambda\to0$, we have $p_*\to1$, so
\[
M_i\longrightarrow
\begin{pmatrix}
5&2\\
2&1
\end{pmatrix}
=A_0,
\qquad i=1,2.
\]
Hence $\log\rho(M_i)\to\log\Lambda_0$. Since
\[
\log\rho(M_i)
\le h_\mu(T_\lambda^{(i)})
\le\log\Lambda_0,
\]
we conclude that
\[
h_\mu(T_\lambda^{(i)})
\longrightarrow\log\Lambda_0
\qquad\text{as }\lambda\to0,
\qquad i=1,2.
\]
\end{proof}

%==============================================================================
%  Appendices
%==============================================================================
\appendix
\section{Katznelson's Diophantine Estimate}
\label{chpt:katznelson-appendix}

Let $B\in\mathrm{M}(d,\mathbb Z)$ with $\det B\neq0$.  We use the notation $E^{cs}$, $E^u$, $\Lambda$,
$p_{cs}$, $r$, and $q$ introduced in
Section~\ref{ch3:auto}, and write $k=k_{cs}+k_u$ according to the
corresponding spectral splitting.

\begin{lemma}[Structural lemma]
\label{structural lemma}
Assume that $E^u\neq\{0\}$.
Then there exists $C_u>0$ such that
\begin{equation}
\label{d-dim:unstable-jordan-lower-bound}
    |B^nv|
    \geq C_u n^{-q}\Lambda^n|v|,
    \qquad v\in E^u,\quad n\geq1.
\end{equation}
\end{lemma}

\begin{proof}
On a Jordan block of size $s$ associated with an eigenvalue $\lambda$ of
modulus greater than one,
\[
    \|B^{-n}\|\leq C n^{s-1}|\lambda|^{-n}.
\]
On the blocks for which $|\lambda|=\Lambda$, the largest polynomial
factor is $n^q$.  Polynomial factors belonging to blocks with
$|\lambda|>\Lambda$ are absorbed into the exponential gap
$|\lambda|/\Lambda>1$.  Hence
\[
    \|(B|_{E^u})^{-n}\|\leq C n^q\Lambda^{-n}.
\]
Applying this estimate to $B^nv$ and taking reciprocals proves
\eqref{d-dim:unstable-jordan-lower-bound}.
\end{proof}

If $B$ is ergodic, then $E^u\neq\{0\}$.  Indeed, otherwise all
eigenvalues would have modulus at most one.  Since $\det B$ is a nonzero
integer, the product of their moduli is at least one, and consequently all
eigenvalues would have modulus one.  Kronecker's theorem would then imply
that they are roots of unity, contradicting ergodicity.

\begin{lemma}[Katznelson's lemma; \cite{Kat71}]
\label{Katznelson}
If $B\in \mathrm{M} (d,\mathbb Z)$ with $\det B\neq0$ is ergodic, then there exists $c_{\mathrm{K}}=c_{\mathrm{K}}(B)>0$ such that
\begin{equation*}
    |k_u|\geq c_{\mathrm{K}}|k|^{-r},
    \qquad 0\neq k\in\mathbb Z^d.
\end{equation*}
\end{lemma}
Neither irreducibility of the characteristic polynomial nor
diagonalizability of $B$ is required.

\section{Invariance of Lebesgue Measure}
\begin{comment}
We record the elementary criterion used in the nonlinear examples.  For
$k=(k_1,\dots,k_d)\in\mathbb Z^d$ and
$z=(z_1,\dots,z_d)\in\mathbb T^d$, recall that
$z^k=z_1^{k_1}\cdots z_d^{k_d}$.
\end{comment}
\begin{lemma}[Fourier criterion for Lebesgue invariance]
\label{lm:Leb meas preserv}
Let $T:\mathbb T^d\to\mathbb T^d$ be a Borel map.  Then $T$ preserves
normalized Lebesgue measure $\mu$ if and only if
\[
   \int_{\mathbb T^d}(Tz)^k\,d\mu(z)=\delta_{k,0}
   \qquad\text{for every }k\in\mathbb Z^d.
\]
\end{lemma}

\begin{proof}
Let $\nu=T_*\mu$.  By the definition of pushforward measure,
\[
   \int_{\mathbb T^d}y^k\,d\nu(y)
   =\int_{\mathbb T^d}(Tz)^k\,d\mu(z).
\]
The Fourier coefficients of normalized Lebesgue measure are
\[\int_{\mathbb T^d}y^k\,d\mu(y)=\delta_{k,0}.\]  Thus the displayed
condition is equivalent to equality of all Fourier coefficients of $\nu$ and
$\mu$.  Since trigonometric polynomials are uniformly dense in
$C(\mathbb T^d)$, this is equivalent to $\nu=\mu$, or $T_*\mu=\mu$.
\end{proof}
\section{Mixing for Higher-Rank Toral Actions}
\label{chpt:higher-rank-actions}

Below we apply the Fourier escape argument from
Section~\ref{section:analytic Fourier framework} to higher-rank actions.  Throughout this
appendix, $|\cdot|=|\cdot|_1$ on both $\mathbb Z^d$ and $\mathbb Z^k$.
Let $T_1,\ldots,T_k$ be commuting toral maps induced by integer matrices
and suppose that each $T_i$ preserves $\mu$.  

Here and below, $\mathbf n\in\mathbb Z^k$ in the automorphism case and
$\mathbf n\in\mathbb Z_+^k$ in the endomorphism case.  Put
\[
 T^{\mathbf n}:=T_1^{n_1}\cdots T_k^{n_k},
 \qquad
 B^{\mathbf n}:=(T^{\mathbf n})^{\mathsf t}.
\]
Then $T^{\mathbf n}$ preserves $\mu$ and, for every
$m\in\mathbb Z^d$,
\[
 (T^{\mathbf n}z)^m=z^{B^{\mathbf n}m}.
\]

We define $\mathcal C_{f,g}(\mathbf n)$ as in
Section~\ref{section:analytic Fourier framework}, with $T^n$ replaced by
$T^{\mathbf n}$.

\begin{lemma}
\label{higher-rank:orbit-growth-implies-mixing}
Suppose that there are
$c_0,\tau>0$ and $s\geq0$ such that
\begin{equation}
\label{higher-rank:abstract-orbit-growth}
 |B^{\mathbf n}m|
 \geq c_0e^{\tau|\mathbf n|}|m|^{-s},
 \qquad 0\neq m\in\mathbb Z^d.
\end{equation}
Then, for every $0<\theta<1$, there are $C_{\theta,d},\beta_\theta>0$
such that
\begin{equation}
\label{higher-rank:abstract-correlation-bound}
 \mathcal C_{f,g}(\mathbf n)
 \leq C_{\theta,d}
 \exp\left(-\beta_\theta
 e^{\frac{\tau}{s+1}|\mathbf n|}\right)
 \|f\|_\theta\|g\|_\theta
\end{equation}
for all $f,g\in\mathcal H_\theta$.
\end{lemma}

\begin{proof}
Set
\[
 R_{\mathbf n}
 :=\left(\frac{c_0}{2}\right)^{1/(s+1)}
 e^{\frac{\tau}{s+1}|\mathbf n|}.
\]
If $0<|m|<R_{\mathbf n}$, then
\eqref{higher-rank:abstract-orbit-growth} implies
\[
 |B^{\mathbf n}m|
 >c_0e^{\tau|\mathbf n|}R_{\mathbf n}^{-s}
 =2R_{\mathbf n}.
\]
Thus, no nonconstant character with frequency below $R_{\mathbf n}$ is
sent to one with frequency below $R_{\mathbf n}$.  Applying
Lemma~\ref{lem:Fourier escape argument} to the single map
$T^{\mathbf n}$ gives
\[
 \mathcal C_{f,g}(\mathbf n)
 \leq C_{\theta,d}\theta^{R_{\mathbf n}/d}
 \|f\|_\theta\|g\|_\theta,
\]
which is \eqref{higher-rank:abstract-correlation-bound} with
\[
 \beta_\theta
 :=\frac{-\log\theta}{d}
 \left(\frac{c_0}{2}\right)^{1/(s+1)}.
\]
\end{proof}

\subsection{Actions by Automorphisms}

Let $T_1,\ldots,T_k\in\mathrm{GL}(d,\mathbb Z)$.  The resulting
$\mathbb Z^k$-action is \emph{totally ergodic} if $T^{\mathbf n}$ is
ergodic for every $0\neq\mathbf n\in\mathbb Z^k$. 

The estimate \eqref{higher-rank:abstract-orbit-growth} with $s=d$ is
proved in \cite[Theorem~2.1]{KK05} and restated in
\cite[Lemma~4.3(a)]{DK10}; see also the box-escape form in
\cite[Lemmas~4.1--4.5]{FKS13}.

\begin{theorem}
Every totally ergodic $\mathbb Z^k$-action by toral automorphisms is
mixing at a double exponential rate on $\mathcal H_\theta$, uniformly in
$\mathbf n$.  More precisely, there exists $\tau>0$ such that, for every
$0<\theta<1$, there are $C_{\theta,d},\beta_\theta>0$ for which
\begin{equation*}
 \mathcal C_{f,g}(\mathbf n)
 \leq C_{\theta,d}
 \exp\left(-\beta_\theta
 e^{\frac{\tau}{d+1}|\mathbf n|}\right)
 \|f\|_\theta\|g\|_\theta
\end{equation*}
for all $f,g\in\mathcal H_\theta$ and $\mathbf n\in\mathbb Z^k$.
\end{theorem}

\begin{proof}
Apply Lemma~\ref{higher-rank:orbit-growth-implies-mixing} with $s=d$.
\end{proof}

\subsection{Actions by Endomorphisms}

Now let $T_1,\ldots,T_k\in\mathrm M(d,\mathbb Z)$ have nonzero
determinants.  Each $T_i$ induces a surjective endomorphism of
$\mathbb T^d$.  The $\mathbb{Z}_+^k$-action is
\emph{totally ergodic} if every nonzero element of the induced $\mathbb Z^k$-action on
its natural invertible extension is ergodic.

The estimate \eqref{higher-rank:abstract-orbit-growth} with some $s>0$ is
proved in \cite[Theorem~4]{L13}.  

\begin{theorem}
Every totally ergodic $\mathbb Z_+^k$-action by nonsingular toral
endomorphisms is mixing at a double exponential rate on
$\mathcal H_\theta$, uniformly in $\mathbf n$.  More precisely, there are
$\tau,s>0$ such that, for every $0<\theta<1$, there are
$C_{\theta,d},\beta_\theta>0$ for which
\begin{equation*}
 \mathcal C_{f,g}(\mathbf n)
 \leq C_{\theta,d}
 \exp\left(-\beta_\theta
 e^{\frac{\tau}{s+1}|\mathbf n|}\right)
 \|f\|_\theta\|g\|_\theta
\end{equation*}
for all $f,g\in\mathcal H_\theta$ and
$\mathbf n\in\mathbb Z_+^k$.
\end{theorem}

\begin{proof}
Apply Lemma~\ref{higher-rank:orbit-growth-implies-mixing}.
\end{proof}
\section{Higher-Dimensional Geometric Series}

The following estimates are used in the Fourier escape argument.

\begin{lemma}
\label{lem:higher-dimensional-geometric-series}
Let $0<\eta<1$ and $R>0$.  Then
\begin{equation}
\label{appendix:full geometric series}
 \sum_{k\in\mathbb Z^d\setminus\{0\}}\eta^{|k|_1}
 =\left(\frac{1+\eta}{1-\eta}\right)^d-1
 =:S_{\eta,d},
\end{equation}
and
\begin{equation*}
 \sum_{|k|_1\ge R}\eta^{|k|_1}
 \le D_{\eta,d}\eta^{R/d},
 \qquad
 D_{\eta,d}:=
 2d\left(\frac{1+\eta}{1-\eta}\right)^d\frac1{1+\eta}.
\end{equation*}
\end{lemma}

\begin{proof}
Separating the coordinates gives
\[
 \sum_{k\in\mathbb Z^d}\eta^{|k|_1}
 =\prod_{s=1}^d\sum_{m\in\mathbb Z}\eta^{|m|}
 =\left(\frac{1+\eta}{1-\eta}\right)^d,
\]
which proves \eqref{appendix:full geometric series}.  Since
$|k|_\infty\ge |k|_1/d$, the set $\{|k|_1\ge R\}$ is contained in the
union of the $2d$ half-spaces $k_s\ge R/d$ and $k_s\le-R/d$.  Therefore,
by symmetry,
\begin{align*}
 \sum_{|k|_1\ge R}\eta^{|k|_1}
 &\le 2d\sum_{k_1\ge R/d}\eta^{|k|_1} \\
 &\le \frac{2d\eta^{R/d}}{1-\eta}
       \left(\frac{1+\eta}{1-\eta}\right)^{d-1}
 =D_{\eta,d}\eta^{R/d}.
\end{align*}
\end{proof}

\section{Mixing Hierarchy}

Below, we show that mixing at a double exponential rate on an analytic space implies superexponential mixing on $C^\infty$. The proof adapts the approximation arguments from
\cite[Lemma~C.1]{Shc26} and \cite[Lemma~2.4]{GS14}.

For an integer $m\geq0$, let 
\[
    \|h\|_{H^m}^2
    :=\sum_{k\in\mathbb Z^d}(1+|k|)^{2m}|h_k|^2
\]
be the Sobolev norms on $\mathbb T^d$.

\begin{lemma}
\label{d-dim:analytic-to-smooth}
Let $T:\mathbb T^d\to\mathbb T^d$ preserve normalized Lebesgue measure.
Suppose that, for some $0<\theta<1$ and some $C,\beta,\gamma>0$,
\begin{equation}
\label{d-dim:abstract-analytic-double-exponential}
    \mathcal C_{f,g}(n)
    \leq
    C\exp\bigl(-\beta e^{\gamma n}\bigr)
    \|f\|_\theta\|g\|_\theta
\end{equation}
for all $f,g\in\mathcal H_\theta$ and all $n\geq0$.
Then, for every $\alpha>0$, there exist an integer $m\geq0$ and a constant
$C_\alpha>0$ such that
\begin{equation}
\label{d-dim:smooth-superexponential}
    \mathcal C_{f,g}(n)
    \leq C_\alpha e^{-\alpha n}\|f\|_{C^m}\|g\|_{C^m}
\end{equation}
for all $f,g\in C^\infty(\mathbb T^d)$ and all $n\geq0$.
One may take any integer $m>\alpha/\gamma$.
\end{lemma}

\begin{proof}
Subtracting the averages from $f$ and $g$ does not change their
correlation.  We may therefore suppose that $f_0=g_0=0$, so that
\[
    \mathcal C_{f,g}(n)
    =\bigl|\langle f\circ T^n,g\rangle_{L^2}\bigr|.
\]
For $N\geq0$, denote the Fourier truncation of $h$ by
\[
    P_Nh:=\sum_{|k|\leq N}h_kz^k.
\]
Put $\rho:=-\log\theta>0$.  Since
$\theta^{-|k|}=e^{\rho|k|}\leq e^{\rho N}$ for $|k|\leq N$,
\begin{equation}
\label{d-dim:truncation-analytic-norm}
    \|P_Nh\|_\theta
    \leq e^{\rho N}\|P_Nh\|_{L^2}
    \leq e^{\rho N}\|h\|_{L^2}.
\end{equation}
On the other hand, for every integer $m\geq0$,
\begin{equation}
\label{d-dim:truncation-smooth-tail}
    \|(I-P_N)h\|_{L^2}
    \leq(1+N)^{-m}\|h\|_{H^m}.
\end{equation}

Choose
\[
    \delta:=\frac{\beta}{4\rho}
\]
and, for all sufficiently large $n$, set
\[
    N_n:=\left\lfloor\delta e^{\gamma n}\right\rfloor.
\]
After increasing the lower threshold for $n$, we have
\[
    \frac{\delta}{2}e^{\gamma n}\leq N_n
    \leq\delta e^{\gamma n}.
\]
Applying \eqref{d-dim:abstract-analytic-double-exponential} to
$P_{N_n}f$ and $P_{N_n}g$, and then using
\eqref{d-dim:truncation-analytic-norm}, gives
\begin{align*}
\bigl|\langle (P_{N_n}f)\circ T^n,P_{N_n}g\rangle_{L^2}\bigr|
&\leq
C\exp\bigl(-\beta e^{\gamma n}+2\rho N_n\bigr)
\|f\|_{L^2}\|g\|_{L^2}\nonumber\\
&\leq
C\exp\left(-\frac{\beta}{2}e^{\gamma n}\right)
\|f\|_{L^2}\|g\|_{L^2}.
\end{align*}

Because $T$ preserves Lebesgue measure, its Koopman operator is an
isometry on $L^2$.  Hence
\begin{align*}
&\left|
\langle f\circ T^n,g\rangle_{L^2}
-\langle(P_{N_n}f)\circ T^n,P_{N_n}g\rangle_{L^2}
\right|\\
&\qquad\leq
\|(I-P_{N_n})f\|_{L^2}\|g\|_{L^2}
+\|P_{N_n}f\|_{L^2}\|(I-P_{N_n})g\|_{L^2}.
\end{align*}
Using \eqref{d-dim:truncation-smooth-tail},
$\|P_{N_n}f\|_{L^2}\leq\|f\|_{L^2}$, and the lower bound on $N_n$, we
obtain
\begin{align}
\mathcal C_{f,g}(n)
&\leq
C\exp\left(-\frac{\beta}{2}e^{\gamma n}\right)
\|f\|_{L^2}\|g\|_{L^2}\nonumber\\
&\quad+
C_{m,\delta}e^{-m\gamma n}
\left(
  \|f\|_{H^m}\|g\|_{L^2}
  +\|f\|_{L^2}\|g\|_{H^m}
\right)
\label{d-dim:analytic-smooth-interpolation}
\end{align}
for all sufficiently large $n$.

Given $\alpha>0$, choose an integer $m>\alpha/\gamma$.  Both terms on the
right-hand side of
\eqref{d-dim:analytic-smooth-interpolation} are then bounded by a
constant times $e^{-\alpha n}$.  Enlarging the constant handles the finitely
many remaining values of $n$.  Finally,
\[
    \|h\|_{H^m}\leq C_{m,d}\|h\|_{C^m}
\]
on $\mathbb T^d$, and centering a function changes its $C^m$ norm by at
most a fixed factor.  This proves
\eqref{d-dim:smooth-superexponential}.
\end{proof}

\begin{corollary}
\label{d-dim:toral-smooth-superexponential}
Every ergodic toral endomorphism covered by
Theorem~\ref{d-dim:ergodic auto is dem} is superexponentially mixing on
$C^\infty(\mathbb T^d)$ in the following quantitative sense: for every
$\alpha>0$, there exist $m\geq0$ and $C_\alpha>0$ for which
\eqref{d-dim:smooth-superexponential} holds.  More precisely, any integer
$m$ satisfying
\[
    m>\frac{\alpha(r+1)}{\log\Lambda}
\]
is sufficient.
\end{corollary}

\begin{proof}
For such an $m$, choose
\[
    \frac{\alpha}{m}<\gamma<
    \frac{\log\Lambda}{r+1}.
\]
Theorem~\ref{d-dim:ergodic auto is dem} and the discussion following it
give \eqref{d-dim:abstract-analytic-double-exponential} for this
$\gamma$, while $m>\alpha/\gamma$.  The conclusion therefore follows from
Lemma~\ref{d-dim:analytic-to-smooth}.
\end{proof}

For comparison, Wang \cite[Theorem~10.1]{Wan25} obtains
quantitative superexponential mixing estimates for ergodic
toral automorphisms as part of a broader theory for
nilmanifold automorphisms.

%==============================================================================
%  Bibliography
%==============================================================================
\bibliographystyle{amsalpha}
\bibliography{Katia}

% biblatex users: replace the two lines above with
%\printbibliography
\end{document}